# Omega-inconsistency in Gödel's formal system: a constructive proof of the Entscheidungsproblem

*Bhupinder Singh Anand*

*If we apply an extension of the Deduction meta-Theorem to Gödel's meta-reasoning of "undecidability", we can conclude that Gödel's formal system of Arithmetic is not omega-consistent. If we then interpret [(Ax)F(x)] as "There is a general, x-independent, routine to establish that F(x) holds for all x", instead of as "F(x) holds for all x", it follows that a constructively interpreted omega-inconsistent system proves Hilbert's Entscheidungsproblem negatively.*

**1.1 Notation**

We generally follow the notation of Gödel [Go31a]. However, we use the notation "(A*x*)", which classically interprets as "for all *x*", to denote Gödel's special symbolism for Generalisation.

We use square brackets to indicate that the expression (including square brackets) only denotes the string[1] named within the brackets. Thus, "[(A*x*)]" is not part of the formal system P, and would be replaced by Gödel's special symbolism for Generalisation in order to obtain the actual string in which it occurs.

Following Gödel's definitions of well-formed formulas[2], we note that juxtaposing the string "[(A*x*)]" and the formula[3] "[*F*(*x*)]" is the formula "[(A*x*)*F*(*x*)]", juxtaposing the symbol "[~]" and the formula "[*F*]" is the formula "[~*F*]", and juxtaposing the symbol "[v]" between the formulas "[*F*]" and "[*G*]" is the formula "[*F*v*G*]".

---

[1] We define a "string" as any concatenation of a finite set of the primitive symbols of the formal system under consideration.

[2] We note that all "well-formed formulas" of P are "strings" of P, but all "strings" of P are not "well-formed formulas" of P.

The numerical functions and relations in the following are defined explicitly by Gödel [Go31a]. The formulas are defined implicitly by his reasoning.

**1.2 Definitions**

We take P to be Gödel's formal system, and define ([Go31a], Theorem VI, p24-25):

(*i*)   $Q(x, y)$ as Gödel's recursive numerical relation $\sim xB(Sb(y\ 19|Z(y)))$.

(*ii*)  $[R(x, y)]$ as a formula that represents $Q(x, y)$ in the formal system P.

(The existence of such a formula follows from Gödel's Theorem VII [Go31a].)

(*iii*) $q$ as the Gödel-number of the formula $[R(x, y)]$ of P.

(*iv*)  $p$ as the Gödel-number of the formula $[(Ax)][R(x, y)]$[4] of P.

(*v*)   $[p]$ as the numeral that represents the natural number $p$ in P.

(*vi*)  $r$ as the Gödel-number of the formula $[R(x, p)]$ of P.

(*vii*) $17Genr$ as the Gödel-number of the formula $[(Ax)R(x, p)]$ of P.

(*viii*) $Neg(17Genr)$ as the Gödel-number of the formula $[\sim(Ax)R(x, p)]$ of P.

---

[3] By "formula", we shall mean a "well-formed formula" as defined by Gödel.

[4] We note that "$[(Ax)][R(x, y)]$" and "$[(Ax)R(x, y)]$" denote the same formula of P.

(*ix*)  *R*(*x*, *y*) as the arithmetical interpretation of the formula [*R*(*x*, *y*)] of P.

(*R*(*x*, *y*) is defined by Gödel's Theorem VII [Go31a], where it is proved instantiationally equivalent to *Q*(*x*, *y*).)

## 1.3 Gödel's Lemmas

In the proof of his Theorem VI, Gödel [Go31a] proves the following lemmas:

**Lemma 1**:   ~*xB*(17*Genr*) => *Bew*(*Sb*(*r* 17|*Z*(*x*)))

**Lemma 2**:   *xB*(17*Genr*) => *Bew*(*Neg Sb*(*r* 17|*Z*(*x*)))

## 1.4 Gödel's meta-Lemmas

He then proves the following meta-lemmas:

**Meta-lemma 1**:   *Bew*(17*Genr*) => *Bew*(*Neg Sb*(*r* 17|*Z*(*n*))) holds for some natural number *n*.

**Meta-lemma 2**:   *Bew*(17*Genr*) => ~*Bew*(17*Genr*) holds, if P is assumed consistent.

**Meta-Lemma 3**:   ~*nB*(17*Genr*) holds for any natural number *n*, if P is assumed consistent.

**Meta-lemma 4**:   *Bew*(*Neg*(17*Genr*)) => ~*Bew*(*Neg*(17*Genr*)) holds, assuming P is omega-consistent.

## 1.5 Gödel's conclusions and consequences

From these Gödel concludes that ~*Bew*(17*Genr*) holds if P is assumed consistent. Hence 17*Gen* is not P-PROVABLE[5], and, ipso facto, [(A$x$)$R$($x$, $p$)] is not P-provable, if P is assumed consistent.

He also concludes that ~*Bew*(*Neg*(17*Genr*)) holds if P is assumed omega-consistent. Hence *Neg*(17*Genr*) is not P-PROVABLE, and, ipso facto, [~(A$x$)$R$($x$, $p$)] is not P-provable, if P is assumed omega-consistent.

## 1.6 P is not omega-consistent

**Meta-Theorem 1**: P is not omega-consistent

**Proof**: Since 17*Genr* is the Gödel-number of the formula [(A$x$)$R$($x$, $p$)]:

(*i*)  If [(A$x$)$R$($x$, $p$)] is P-provable, then $nB$(17*Genr*) holds for some natural number $n$.[6]

(*ii*)  Hence, by *Lemma 2*:

---

[5] The web-version [Go31b] of Gödel's paper uses *italics* instead of CAPITALS to refer to meta-mathematical concepts in assertions where the formulas of P are referred to by their Gödel-numbers.

[6] In a companion paper [An02], where I review Gödel's and Rosser's non-formal meta-reasoning of undecidability, I argue that this "semantic" meta-equivalence is equivalent to the "non-semantic" meta-assertion:

   [(A$x$)$R$($x$, $p$)] => [*Bew*(*Neg Sb*($r$ 17|$Z$($n$)))] is P-provable for some numeral [$n$],

where [*Bew*(*Neg Sb*($r$ 17|$Z$($n$)))] is the formula of P whose interpretation is the proposition obtained when we substitute a given natural number $n$ for the variable $x$ in the recursive relation *Bew*(*Neg Sb*($r$ 17|$Z$($x$))). I argue that §1.6(*ii*) and (*iii*) can similarly be expressed as "non-semantic" P-provable meta-assertions.

If [(A*x*)*R*(*x*, *p*)] is P-provable, then *Bew*(*Neg Sb*(*r* 17|*Z*(*n*))) holds for some natural number *n*.

(*iii*) Since *Neg Sb*(*r* 17|*Z*(*n*)) is the Gödel-number of [~*R*(*n*, *p*)], we have that:

If *Bew*(*Neg Sb*(*r* 17|*Z*(*n*))) holds for some natural number *n*, then [~*R*(*n*, *p*)] is P-provable for some numeral [*n*].

(*iv*) We thus have the meta-inferences:

If [(A*x*)*R*(*x*, *p*)] is P-provable, then [~*R*(*n*, *p*)] is P-provable for some numeral [*n*]

If [(A*x*)*R*(*x*, *p*)] is P-provable, then [(E*x*)~*R*(*x*, *p*)] is P-provable[7]

If [(A*x*)*R*(*x*, *p*)] is P-provable, then [~(A*x*)*R*(*x*, *p*)] is P-provable

(*v*) We now appeal to an extension of the Deduction Theorem (see Appendix 1), and conclude that:

[(A*x*)*R*(*x*, *p*)] => [~(A*x*)*R*(*x*, *p*)] is P-provable.

(*vi*) By the logical axioms of P, it follows that:

[~(A*x*)*R*(*x*, *p*)] is P-provable

(*vii*) Hence [(A*x*)*R*(*x*, *p*)] is not P-undecidable.

(*viii*) Now, from *Meta-lemma 3* we have:

[*R*(*n*, *p*)] is P-provable for some numeral [*n*].

---

[7] We use "(E*x*)" to denote the string "~(A*x*)~".

(*ix*)   It follows from *(vi)* and *(viii)* that P is not omega-consistent.

**1.7 Conclusion**

We conclude that, if we admit meta-mathematical arguments of provability, then P is not omega-consistent (we note an interesting interpretation of this in Appendix 2).

**Appendix 1: An "extended" Deduction Theorem**

In §1.6(v), we appeal to the following argument for an "extended" interpretation of the Deduction Theorem.

**Deduction Theorem**: From "If [*A*] is P-provable, then [*B*] is P-provable" we may conclude "([*A*] => [*B*]) is P-provable", where [*A*], [*B*] are propositions in P.

**Proof**: The meta-deduction "If [*A*] is P-provable, then [*B*] is P-provable", where [*A*], [*B*] are propositions in P, implies that there is some finite sequence of P-formulas, [$B_1$], [$B_2$],..., [$B_n$], such that [$B_1$] is [*A*], and each |$B_i$| is (*1*) either [*A*], or (*2*) an axiom of P, or (*3*) an immediate consequence of some formulas in the set consisting of [*A*], the axioms of P, and the preceding [$B_i$]s. Using induction on "*i*", we may thus conclude that ([*A*] => [*B*]) is P-provable by the usual reasoning ([Me64], p61, proposition 2.4).

We highlight the distinction, since it can be argued that the proof of the standard Deduction Theorem requires that, in order to conclude that ([*A*] => [*B*]) is P-provable, the sequence [$B_1$], [$B_2$],..., [$B_n$] must be explicitly expressed, and not simply assumed as implicitly expressible.

Now, by Gödel's reasoning, *Bew*(*n*) may hold for a class *D* of natural numbers even when an explicit PROOF for each *n* in the class is not known. In other words, we may have no known general method for constructively determining the PROOF of any given *n* of *D*. This simply reflects the existential nature of Gödel's definition *Bew*(*x*) <=> (*Ey*)*yBx*.

Hence, an extended interpretation of the Deduction Theorem becomes necessary, unless we deny validity to meta-mathematical reasoning for establishing that a formula, or a class of formulas, of P is provable. However, this would then deny validity to Gödel's meta-mathematical Lemmas, from which he concludes his Theorem VI.

Appeal to an extended Deduction Theorem is also explicit in Mendelson's version ([Me64], p146, proposition 3.32) of Rosser's argument[8] for the construction of undecidable propositions in a consistent P.

## Appendix 2: A constructive interpretation of "(Ax)F(x)" and a solution to Hilbert's Entscheidungsproblem for P

The omega-inconsistency of P has an interesting, constructive, interpretation that yields a negative proof of Hilbert's Entscheidungsproblem for P.

*The classical Platonist interpretation*

It follows from §1.6(*viii*) that the interpretation *R*(*n*, *p*) of [*R*(*n*, *p*)] holds for any natural number *n*. It follows from §1.6(*vi*) that the interpretation ~(A*x*)*R*(*x*, *p*) of [~(A*x*)*R*(*x*, *p*)] also holds. If we ignore intuitionist objections, and interpret the latter as the Platonist meta-assertion:

   There is some natural number *n* for which *R*(*n*, *p*) does not hold,

---

[8] However, Mendelson's version appears to apply the extended Deduction Theorem invalidly. Assuming that *r* is a PROOF of a given FORMULA *n*, it seems to make the invalid assumption, "[*r*]=<*x* is P-provable", as the premise in the application of the extended Deduction Theorem.

(We use the notation "=<" to to denote the symbol that interprets as "equal to or less than".)

Then **Meta-theorem 1** has the uncomfortable consequence that the standard interpretation of P is inconsistent.

*A constructive interpretation of Generalisation*

However, we can also consider the constructive, and intuitionistically unobjectionable, interpretation $(Ax)F(x)$ as the meta-assertion:

There is a general, $x$-independent, routine to establish that $F(x)$ holds for all $x$.

In other words, we take the standard interpretation of $[(Ax)F(x)]$ as the assertion that we can always construct a Turing machine T, independent of $n$, which can decide that $F(n)$ holds for any given natural number $n$.

We note that Generalisation ([Me64], p57) would then interpret as a constructive, and intuitionistically unobjectionable, Rule of Inference.

*Hilbert's Entscheidungsproblem*

So, if we interpret P constructively, we would conclude from §1.6(*vi*) and §1.6(*viii*) that whereas, given any natural number $n$, we can always find some $n$-dependent method to establish that $R(n, p)$ holds, we cannot find a general $n$-independent method to establish that $R(n, p)$ holds for any, or all $n$.

In other words, given any natural number $n$, we can always construct a Turing machine $T(n)$, that depends on $n$, which will decide whether the interpretation $R(n, p)$ of $[R(n, p)]$ holds or not. However, we cannot construct a Turing machine T that is independent of $n$, and which will decide, for any given $n$ as input, whether $R(n, p)$ holds or not.

Thus the omega-inconsistency of P can be seen as a constructive, and intuitionistically unobjectionable, negative proof of Hilbert's Entscheidungsproblem under a constructive standard interpretation for P.

# References

[An01]   Anand, Bhupinder Singh. 2001. *Beyond Gödel: Simply constructive systems of first order Peano's Arithmetic that do not yield undecidable propositions by Gödel's reasoning*. Alix_Comsi, Mumbai (*Unpublished*).

<*Web page*: http://alixcomsi.com/Constructivity_abstracts.htm>

[An02]   Anand, Bhupinder Singh. 2002. *Reviewing Gödel's and Rosser's non-formal meta-reasoning of undecidability*. Alix_Comsi, Mumbai (*Unpublished*).

<*Web page*: http://www.alixcomsi.com/Constructivity_consider.htm>

[Go31a]  Gödel, Kurt. 1931. *On formally undecidable propositions of Principia Mathematica and related systems I*. In M. Davis (ed.). 1965. The Undecidable. Raven Press, New York.

[Go31b]  Gödel, Kurt. 1931. *On formally undecidable propositions of Principia Mathematica and related systems I*.
<Web-page: http://www.ddc.net/ygg/etext/godel/godel3.htm>

[Me64]   Mendelson, Elliott. 1964. Introduction to Mathematical Logic. Van Norstrand, Princeton.

(*Acknowledgement: I thank Dr. Damjan Bojadziev, Department of Intelligent Systems, Jozef Stefan Institute, Ljubljana, Slovenia for his encouragement and constructive suggestions regarding the presentation of this paper.*)

(*Updated: Friday 9$^{th}$ May 2003 8:46:47 AM by re@alixcomsi.com*)